\documentclass[12pt]{article}
\usepackage{amssymb,amsmath,cite}
\usepackage[dvips]{graphicx,color}
\usepackage{psfrag,subfigure}
\title{Analytical Solution to Improper Integral of Divergent Power Functions Using The Riemann Zeta Function}
\author{Farhad Aghili\footnote{email: faghili@encs.concordia.ca $\;$ or $\;$ farhad.aghili@gmail.com} and Siamak Tafazoli\footnote{email: s$_{-}$tafazo@encs.concordia.ca} }
\date{}

\begin{document}

\maketitle

\begin{abstract}
This paper presents an analytical  closed-form solution to improper integral $\mu(r)=\int_0^{\infty} x^r dx$, where $r \geq 0$. The solution technique is based on splitting the improper integral into an infinite sum of definite integrals with successive integer limits. The exact solution of every definite integral is obtained by making use of the binomial polynomial expansion, which then  allows expression of the entire summation equivalently  in terms of a weighted sum of Riemann zeta functions. It turns out that the solution fundamentally depends on whether or not $r$ is an integer. If $r$ is a non-negative integer, then the solution is manifested in a finite series of weighted Bernoulli numbers, which is then drastically simplified to a second order rational function  $\mu(r)=(-1)^{r+1}/(r+1)(r+2)$. This is achieved by taking advantage of the relationships between Bernoulli numbers and binomial coefficients. On the other hand,  if $r$ is a non-integer real-valued number, then we prove $\mu(r)=0$ by the virtue of  the elegant relationships between zeta  and gamma functions and their properties.      
\end{abstract}

\section{Introduction}

One of the most  enigmatic things in arithmetic is that sum of an infinite divergent series often results in a finite number. For instance, it has long been established that the summation of all of the Natural numbers is equal to a finite  negative rational number, i.e., 
\begin{equation} \notag
\sum_{n=1}^{\infty} n = -\frac{1}{12}
\end{equation}
Bernard Riemann generalized infinite sum of series over the complex plane    
via the Riemann zeta function 
\begin{equation} \notag
\zeta(s) = \sum_{n=1}^{\infty} \frac{1}{n^s}.
\end{equation}
The zeta function is analytically continued over the entire complex $s$ plane except for $s=1$. 
In spite of being a peculiar and counter-intuitive abstraction in mathematics, divergent series appear also in physics. For instance, it stunningly gives correct results for the Casimir effect, which is the existence of an attractive force between two parallel conducting plates in the vacuum~\cite{Dowling-1989,Schumayer-Hutcheiston-2011}.

Improper integrals of power functions  
\begin{equation} \label{eq:Ir_def}
\mu(r)=\int_0^{\infty} x^r dx 
\end{equation}
are the reminiscent  of the infinite sum of divergent series. These integrals are important in analysis and they also appear in modern physics particularly in calculations of quantum field theory.  The improper integrals cannot be normally computed using the fundamental theorem of calculus and hence the solution requires a novel approach. In this paper, we seek an explicit solution to the improper integrals of power functions by making use of Riemann zeta function and its elegant relationships to the Bernoulli numbers and the gamma function. First, the improper function is split into an infinite sum of definite integrals of which the limits are successive integers, i.e., the difference between the limits of every definite integral is 1. Subsequently, we write the improper function, i.e., the mu function, as an infinite series  by making use of the familiar binomial expansion. It then becomes clear that if the power argument is an integer, the resultant infinite series is equivalent to  a finite series of weighted zeta functions or finite weighted sum of some Bernoulli numbers. Finally, taking advantage of the relationship between  Bernoulli numbers and Binomial coefficient, we derive the analytical solution to be $\mu(r)=(-1)^{r+1}/(r+1)(r+2)$, which is surprisingly simple. However, if the power argument is not an integer, then we show that the improper integral is tantamount to an infinite weighted sum of zeta functions. In the  latter case, we prove that the improper integral is identically zero regardless of the power argument.   

\section{Improper integral of power functions}
Let us define the mu function as the improper integral of power functions in terms of argument $r$  
\begin{equation} \label{eq:Ir_def}
\mu(r)=\int_0^{\infty} x^r dx 
\end{equation}
The mu function can be interpreted as the natural extension of the zeta function where the discrete summation is replaced by a continuous integral. The above improper integral can be equivalently written as the following infinite series by splitting the range of integration limits into successive integer numbers. That is   
\begin{subequations}
\begin{equation} \label{eq:sr=sumDela}
\mu(r)=\sum_{n=1}^{\infty} \Delta(r,n)
\end{equation} 
where 
\begin{equation} \label{eq:Delta_r}
\Delta(r,n)=\int_{n-1}^n x^r dx =\frac{1}{r+1}\Big( n^{r+1} - (n-1)^{r+1}  \Big)
\end{equation}
\end{subequations}
is the determinate integral of the polynomial function over limits $n-1$ and $n$. In the following sections, we will show that the divergent series in the RHS of \eqref{eq:sr=sumDela} converge to a finite solution and that there two set of distinct solutions associated with the cases of the argument $r$ being a positive integer and a non-integer real number.   

\subsection{Domain of natural numbers}

Suppose the domain of mu function \eqref{eq:Ir_def} is the natural numbers, i.e., $r \in \mathbb{N}$, where  $\mathbb{N}=\{0, 1, 2, \cdots \} $. Then, we can write 
the binomial  polynomial expansion of $(n-1)^{r+1}$  by
\begin{subequations}
\begin{equation} \label{eq:n-1_binomial}
(n-1)^{r+1} = \sum_{k=0}^{r+1} C_{r+1}^k  (-1)^{r+1-k}n^k,
\end{equation}
where
\begin{equation} \label{binom}
C_m^l = \frac{m !}{l !(m-l) !}
\end{equation}
\end{subequations}
is  the binomial  coefficient.  Substituting  \eqref{eq:n-1_binomial} into \eqref{eq:Delta_r} yields
\begin{equation} \label{eq:Delta_r2}
\Delta(r,n)=\frac{1}{r+1}\sum_{k=0}^{r} C_{r+1}^k (-1)^{r-k}n^k
\end{equation}
Thus \eqref{eq:sr=sumDela} can be written as
\begin{align} \notag 
\mu(r) & =\frac{1}{r+1}\sum_{n=1}^{\infty} \sum_{k=0}^{r} C_{r+1}^k (-1)^{r-k}n^k \\ \label{eq:Delta_r2}
& = \frac{1}{r+1}  \Big( \sum_{k=0}^{r}(-1)^{r-k} C_{r+1}^k  \Big) \sum_{n=1}^{\infty} n^k  \qquad \forall r \in \mathbb{N} 
\end{align}
The last summation in RHS of \eqref{eq:Delta_r2} can be formally written  
in terms of  the Reimann zeta function $\zeta(\cdot)$ as
\begin{equation} \label{eq:zeta}
\zeta(-k) = \sum_{n=1}^{\infty} n^{k}
\end{equation}
and therefore we arrive at
\begin{equation} \label{eq:sr}
\mu(r)=\frac{1}{r+1}\sum_{k=0}^{r} C_{r+1}^k (-1)^{r-k} \zeta(-k)   \qquad \forall r \in \mathbb{N}
\end{equation}
For positive integer $k\geq 0$, the zeta function is related to Bernoulli numbers by
\begin{equation} \label{eq:zeta_B}
\zeta(-k) = (-1)^k \frac{B_{k+1}}{k+1}
\end{equation}
Moreover, from definition \eqref{binom}, one can verify that the successive binomial coefficients hold the following useful identity   
\begin{equation} \label{eq:2Cs}
C_{r+1}^k  = \frac{k+1}{r+2} C_{r+2}^{k+1} 
\end{equation}
Finally, upon substituting the relevant terms from \eqref{eq:zeta_B} and \eqref{eq:2Cs} into \eqref{eq:sr}, one can equivalently rewrite the latter equation in the following simple form 
\begin{subequations}
\begin{align} \label{eq:srs}
\mu(r)&=\frac{(-1)^{r}}{(r+1)(r+2)}\sum_{k=0}^{r} C_{r+2}^{k+1} B_{k+1}  \qquad \forall r \in \mathbb{N} \\ \label{eq:Sr2}
&= \frac{(-1)^{r}}{(r+1)(r+2)}\sum_{k=1}^{r+1} C_{r+2}^k B_{k}.
\end{align}
\end{subequations}
On the other hand, the bernoulli numbers satisfy the following property 
\begin{equation} \notag
\sum_{k=0}^{m-1} C_m^k B_k =0
\end{equation}
Knowing that $B_0=1$, one can equivalently write  the above equation as 
\begin{equation}  \label{eq:sumCB}
\sum_{k=1}^{r+1} C_{r+2}^k B_k =-1
\end{equation}
Finally upon substituting  \eqref{eq:sumCB} in \eqref{eq:Sr2}, we arrive at the explicit expression of the mu function  
\begin{equation}  \label{eq:Sr3}
\mu(r)=\frac{(-1)^{r+1}}{(r+1)(r+2)}  \qquad \forall r \in \mathbb{N}
\end{equation}

\subsection{Domain of non-integer real numbers}

\begin{figure}
\psfrag{lambda}[c][c][.9]{$\lambda$}
\psfrag{r}[c][c][.9]{$r$} 
\centering{\includegraphics[clip,width=8.5cm]{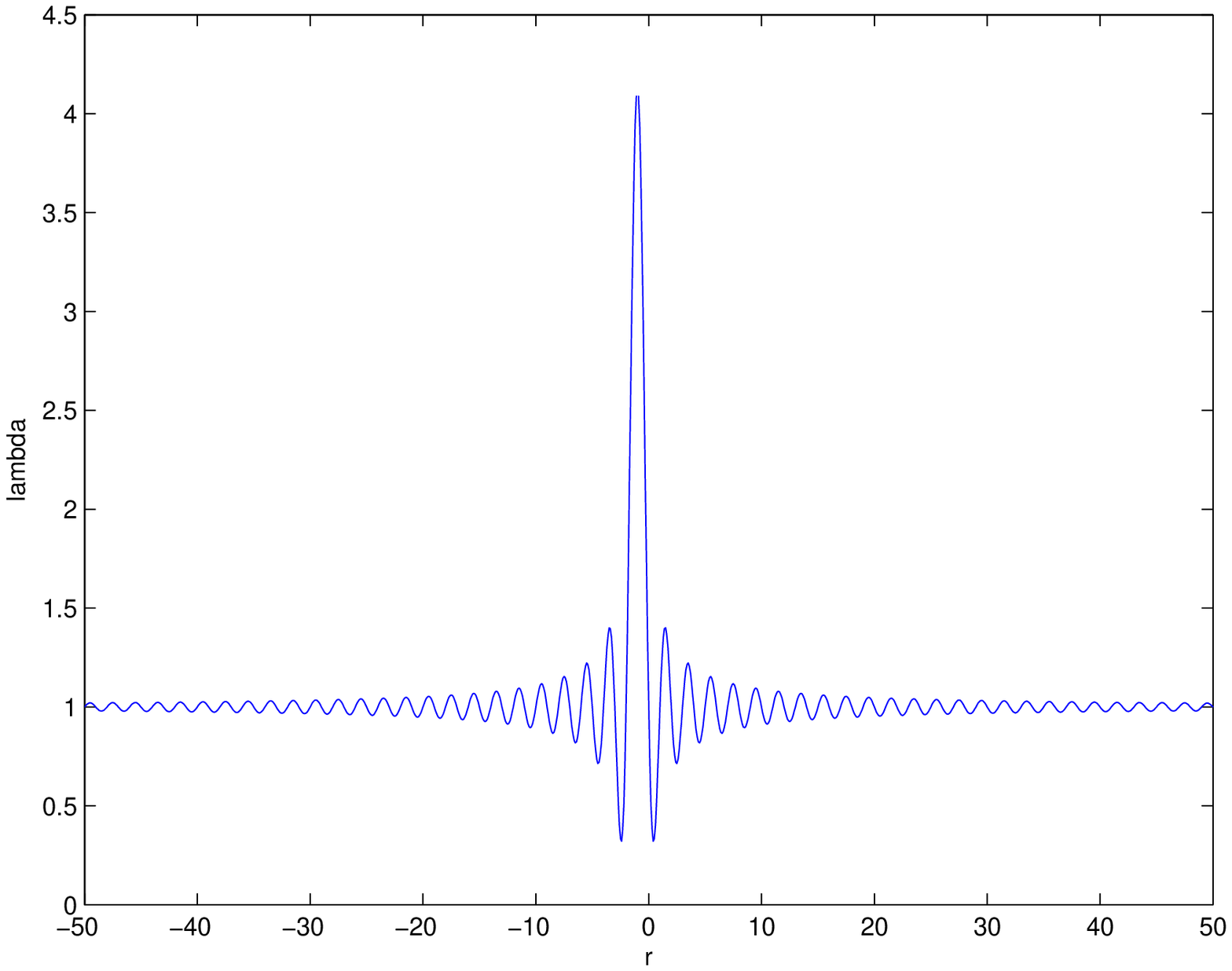}}
\caption{The $\lambda$ function.} \label{fig:lambda}
\end{figure}

Now, consider  the domain of the mu function be real numbers excluding all integers, i.e., $r\in \mathbb{R}-\mathbb{Z}$, where  $\mathbb{Z}=\{\cdots, -2, -1, 0 , 1, 2, \cdots \}$. Writing the Maclaurin expansions about $n=0$, we get
\begin{equation} \label{eq:n_real}
(n-1)^r = \sum_{k=0}^{\infty} C_r^k (-1)^k n^{r-k},    \qquad \forall r\in \mathbb{R}-\mathbb{Z} 
\end{equation}
where the Binomial coefficients in terms of real numbers can be obtained from the gamma function
\begin{equation} \label{eq:Binom_gamma}
C_r^k=\frac{\Gamma(r+1)}{k !  \Gamma(r-k+1) } 
\end{equation}
Notice that the gamma function is not finite at zero and negative integers. Nevertheless the argument of the gamma function in equation \eqref{eq:Binom_gamma}, i.e.,  $r-k+1$, can not be any integer because   $r \in \mathbb{R}-\mathbb{Z} $ while $k \in \mathbb{N}$, and hence the gamma function  is always well defined. Knowing  that 
\begin{equation}
C_{r+1}^0 \equiv 1 \qquad \forall r \in \mathbb{R}-\mathbb{Z} 
\end{equation}
one can write the delta  function  \eqref{eq:Delta_r} in  the following form
\begin{subequations}
\begin{align} \label{eq:Delta_real}
\Delta(r,n)&= \frac{1}{n+1} \Big( n^{r+1} - (n-1)^{r+1} \Big)   \quad \qquad \forall r \in \mathbb{R}-\mathbb{Z} \\
&= \frac{1}{n+1} \Big( n^{r+1} - \sum_{k=0}^{\infty} C_{r+1}^k (-1)^k n^{r+1-k} \Big)\\
&= \frac{1}{n+1} \Big( n^{r+1} -  C_{r+1}^0 n^{r+1}  +\sum_{k=1}^{\infty} C_{r+1}^k (-1)^{k+1} n^{r+1-k} \Big)\\
&= \frac{1}{r+1} \sum_{k=1}^{\infty} C_{r+1}^k (-1)^{k+1} n^{r+1-k} \\
&= \frac{\Gamma(r+2)}{r+1}  \sum_{k=1}^{\infty}\frac{(-1)^{k+1} }{k ! \Gamma(r+2-k)}n^{r+1-k} \\ \label{eq:Delta2}
&= \Gamma(r+1) \sum_{k=1}^{\infty}\frac{(-1)^{k+1} }{k ! \Gamma(r+2-k)}n^{r+1-k} 
\end{align} 
\end{subequations}
where \eqref{eq:Delta2} is obtained using the following property of gamma function
\begin{equation} \label{eq:Gamma_alpha+1}
\Gamma(\alpha +1) = \alpha \Gamma(\alpha)
\end{equation}
Then by virtue of \eqref{eq:sr=sumDela}, the mu function becomes
\begin{subequations}
\begin{align} \notag
\mu(r)&=\Gamma(r+1) \sum_{n=1}^{\infty} \sum_{k=1}^{\infty} \frac{(-1)^{k+1} }{k ! \Gamma(r+2-k)}n^{r+1-k} \\ \label{eq:Ir_sum}
&= \Gamma(r+1) \sum_{k=1}^{\infty} \frac{(-1)^{k+1} \zeta(k-r-1)}{k ! \Gamma(r+2-k)} \qquad \forall r\in \mathbb{R}-\mathbb{Z} 
\end{align}
\end{subequations}
The zeta function can be also expressed in terms of the integral 
\begin{equation} \label{eq:zeta_gamma}
\zeta(s)=\frac{1}{\Gamma(s)} \int_0^{\infty} \frac{t^{s-1}}{e^t -1} dt
\end{equation}
Thus, using \eqref{eq:zeta_gamma} in \eqref{eq:Ir_sum}, we arrive at
\begin{equation} \label{eq:Ir_int}
\mu(r) = \Gamma(r+1)\sum_{k=1}^{\infty} \frac{(-1)^{k+1} }{k ! \Gamma(r-k+2) \Gamma(-r+k-1)} \int_0^{\infty} \frac{t^{-r+k-2}}{e^t -1} dt\end{equation}
On the other hand, the Euler's reflection formula is the fundamental equation of the gamma function \cite{Whittaker-Watson-1920}
\begin{equation} \label{eq:Gamma(1-a)(a)}
\Gamma(1-\alpha) \Gamma(\alpha) = \frac{\pi}{\sin(\pi \alpha)}
\end{equation}
which leads to the following identity
\begin{align} \notag
\frac{1}{\Gamma(r-k+2) \Gamma(-r+k-1)} &= \frac{\sin(\pi (-r+k-1))}{\pi} \\ \label{eq:gamma_sin}
&= \frac{(-1)^k}{\pi}  \sin(\pi r)   \qquad \forall r\in \mathbb{R}-\mathbb{Z} 
\end{align}
Finally, substituting identity  \eqref{eq:gamma_sin} into \eqref{eq:Ir_int} and rearranging the resultant equation, we arrive at
\begin{subequations}
\begin{align}\label{eq:Ir} 
\mu(r)&=- \frac{\sin(\pi r) }{\pi}\Gamma(r+1) \sum_{k=1}^{\infty} \frac{1}{k !} \int_0^{\infty} \frac{t^{k-r-2}}{e^t -1} dt\\ \label{eq:I(r)_int_sum} 
&=-\frac{\sin(\pi r) }{\pi}\Gamma(r+1) \int_0^{\infty} \Big( \sum_{k=1}^{\infty} \frac{t^k}{k !} \Big) \frac{t^{-r-2}}{e^t -1} dt\\
&=-\frac{\sin(\pi r) }{\pi}\Gamma(r+1) \int_0^{\infty}  ( e^t -1 ) \frac{t^{-r-2}}{e^t -1} dt \\
&=-\frac{\sin(\pi r) }{\pi}\Gamma(r+1) \int_0^{\infty}  t^{-r-2} dt \\ \label{eq:I(-r-2)}
&= -\frac{\sin(\pi r) }{\pi}\Gamma(r+1) \mu(-r-2)     \qquad \forall r\in \mathbb{R}-\mathbb{Z} 
\end{align}
\end{subequations}
Note that we justify the interchange of the sum and the integral in \eqref{eq:I(r)_int_sum} based on  uniform convergence assumption of the sum, which will be relaxed at the end where  $\mu(r)$ is proved to be bounded. Now, considering the following change of variable $-r-2$ to replace $r$ in \eqref{eq:I(-r-2)}, one can readily find receptacle of relationship  \eqref{eq:I(-r-2)} as follow
\begin{equation} \label{eq:invIr}
\mu(-r-2)=\frac{\sin(\pi r) }{\pi}\Gamma(-r-1)  \mu(r) 
\end{equation}
It follows by substituting the expression of $\mu(-r-2)$ from \eqref{eq:invIr} into \eqref{eq:I(-r-2)}, i.e.,
\begin{subequations}
\begin{align}
\mu(r)&= - \frac{\sin^2(\pi r)}{\pi^2} \Gamma(r+1) \Gamma (-r-1) \mu(r) \\ \label{eq:Gamma(r+2)}
&= - \frac{\sin^2(\pi r)}{\pi^2 (r+1)} \Gamma(r+2) \Gamma (-r-1) \mu(r) \\ \label{eq:GammaSin}
&=- \frac{\sin^2(\pi r)}{\pi^2 (r+1)} \frac{\pi}{\sin(\pi(- r -1))} \mu(r) \\ \label{eq:twoI(r)}
&= \frac{\sin(\pi r)}{\pi (r+1)} \mu(r)     \qquad \forall r\in \mathbb{R}-\mathbb{Z}, 
\end{align}
\end{subequations}
in which \eqref{eq:Gamma(r+2)} and \eqref{eq:GammaSin} are inferred from identities \eqref{eq:Gamma_alpha+1} and  \eqref{eq:Gamma(1-a)(a)}, respectively. Finally, equation \eqref{eq:twoI(r)} can  be equivalently rewritten in the following forms
\begin{equation} \label{eq:lambad_I}
\lambda(r)  \mu(r) =0   \qquad \forall r\in \mathbb{R}-\mathbb{Z} 
\end{equation}
where 
\begin{equation}
\lambda(r) = 1- \frac{\sin(\pi r)}{\pi(r+1)}.  
\end{equation}
Since $\lambda(r) \neq 0~\forall r\in\mathbb{R}$  (see Fig.~\ref{fig:lambda}),  the only possibility for identity   \eqref{eq:lambad_I} to hold is that
\begin{equation}
\mu(r)=0 \qquad \forall r\in \mathbb{R}-\mathbb{Z},
\end{equation} 
which completes the proof. 

\section{Remarks on the mu function}
One can verify  that 
\begin{equation} 
\sum_{k=0}^{\infty} \mu(k) = \sum_{k=0}^{\infty}  \frac{(-1)^{k+1}}{(k+1)(k+2)}=1 - \log(4)=\log\big(\frac{e}{4} \big).
\end{equation}
Since $\mu(r)=0$ if $r$ takes any non-integer real number, we can conclude 
\begin{equation}
\int_{0}^{\infty} \int_{0}^{\infty} x^r dx dr = \log\big(\frac{e}{4} \big).
\end{equation}

The mu function can be also used to compute the improper integral of a divergent function $\int_0^{\infty} f(x) dx$ using the polynomial series expansion  of the function. For instance, one can shown that
\begin{equation} \notag
\int_0^{\infty} e^x dx = - \frac{1}{e}
\end{equation}

\section{Conclusions}
We derived an explicit solution for the improper integrals of power functions, $\mu(r)=\int_0^{\infty} x^r dx$. Using the binomial expansion, we transcribed the improper integral as an finite or infinite series of weighted zeta functions depending on whether or not the power argument is an integer. By taking advantages of the elegant relationships among zeta function, Bernoulli numbers, binomial coefficient, and the gamma function, the solution of the improper integral has been found to be $\mu(r)=(-1)^{r+1}/(r+1)(r+2)$ if the argument $r$ was an integer, and $\mu(r)=0$ otherwise.   
To this end, it is worth mentioning that although the bounded solution to the the improper integral of divergent power functions is counter-intuitive,  such things need to be understood in the context of the Riemann zeta function pertaining to infinite  divergent series.

\bibliographystyle{IEEEtran}

\end{document}